\documentclass[a4paper]{article}

\usepackage[english]{babel}
\usepackage[utf8x]{inputenc}
\usepackage{amsmath}
\usepackage{amssymb}
\usepackage{amsthm}
\usepackage{graphicx}
\usepackage[colorinlistoftodos]{todonotes}
\usepackage{url}
\usepackage[backref=page]{hyperref}
\usepackage{tensor}
\usepackage[export]{adjustbox}

\title{Spectrally grown graphs}
\author{Mats-Erik Pistol\\ Solid State Physics and NanoLund\\ Box 118, Lund University, S-221 00 Lund SWEDEN\\ \\
Pavel Kurasov\\Department of Mathematics\\ Stockholm University, S-106 91 Stockholm SWEDEN}

\begin{document}
\maketitle

\begin{abstract}
Quantum graphs have attracted attention from mathematicians for some time. A quantum graph is defined by having a Laplacian on each edge of a metric graph and imposing boundary conditions at the vertices to get an eigenvalue problem. A problem studying such quantum graphs is that the spectrum is timeconsuming to compute by hand and the inverse problem of finding a quantum graph having a specified spectrum is difficult. We solve the forward problem, to find the eigenvalues, using a previously developed computer program. We obtain {\bf all} eigenvalues analytically for not too big graphs that have rationally dependent edges. We solve the inverse problem using "spectrally grown graphs". \\  The spectrally grown graphs are evolved from a starting (parent) graph such that the child graphs have eigenvalues are close to some criterion. Our experiments show that the method works and we can usually find graphs having spectra which are numerically close to a prescribed spectrum. There are naturally exceptions, such as if no graph has the prescribed spectrum. The selection criteria (goals) strongly influence the shape of the evolved graphs. Our experiments allows us to make new conjectures concerning the spectra of quantum graphs. We open-source our software at https://github.com/meapistol/Spectra-of-graphs.

\end{abstract}


\section{Introduction}

Quantum graphs have been studied for several decades. A quantum graph \cite{kuchment2004quantum, berkolaiko2013introduction} is a metric graph, $\Gamma$, such that 
with each edge we have associated an ordinary differential equation. Introducing in addition suitable matching conditions at the vertices make the operator self-adjoint.
In the current paper we reduce our analysis to the simplest case of the Laplace differential operator defined on functions satisfying so-called standard vertex conditions (see Eqn. \ref{sc} below).
Then the spectrum of the operator is often refereed to as the spectrum of the graph $\Gamma. $ It is important to 
remember that  there are graphs that are not isomorphic but still have the same spectrum \cite{gutkin2001can}. 
There is a wide literature describing the relations between the topological and metric properties of a finite graph $ \Gamma $ and its spectrum \cite{ berkolaiko2013introduction,kennedy2015spectral}.
For example the spectrum determines the total length of the metric graph and its Euler characteristic \cite{kurasov2008graph,kurasov2008schrodinger}.
Kurasov and Sarnak have recently been shown that exotic crystalline measures can be found using investigations of the additive structure of the spectrum of quantum graphs \cite{kurasov2020stable}

It is not particularly easy to obtain the eigenvalues of a graph which hampers progress. We have previously developed a computer program that gives the eigenvalues of a finite graph \cite{Pistol2016} which was used to find a large set of isospectral graphs \cite{pistol2021generating}. For graphs with pairwise rationally dependent edges (that is, the ratio of the lengths of any two edges is a rational number) we get \emph{all} eigenvalues analytically if the graph is not too big.
For graphs where some edges are not rationally dependent we can only find approximate eigenvalues using a graphical method.
 
Having such a program allows us not only to study spectra of graphs. It also allows us to get a handle of the inverse problem, finding graphs having a prescribed spectrum.
One of the main contributions of our paper is formulation of a new type of evolution that 
we call {\bf spectrum driven evolution}. The spectrum (in our case the eigenvalues) provide an important global characteristic of the system. We find it interesting to define a certain evolution such that the eigenvalues tend to a prescribed set of eigenvalues as the graphs evolve. Other goals can also be given, such as maximizing the first or second spectral gaps. We find that it is indeed possible to evolve graphs in this way such that the goal is closely attained.
%
%
We have found that experimenting with spectral growth of graphs quickly leads to interesting conjectures. One such, concerning the spectral gaps of the first eigenvalues of dumbbell graphs is given is given.

\section{Laplacians on graphs and their spectra}

In this section we briefly define the main object of our studies -- the differential Laplace operator on a metric graph $ \Gamma. $
Consider any finite compact metric graph $ \Gamma $ formed by joining together $ N $ edges $ E_n, \; n=1,2, \dots N $ at $ M $
vertices $ V_m, m=1,2, \dots, M. $ Each edge $ E_n $ is identified with a compact interval $ [x_{2n-1}, x_{2n}] $ on the real line.
Then the (standard) Laplace operator $ L = - \frac{d^2}{dx^2} $ is defined on the functions satisfying standard vertex conditions:
\begin{equation} \label{sc}
\left\{
\begin{array}{ll}
\displaystyle  f(x_i)=f(x_j), & x_i, x_j \in V_m, \\[3mm]
\displaystyle \sum_{x_i{ \in V_m} }\partial_{n}f(x_i)=0. &
\end{array} \right. 
\end{equation}
at each vertex $ V_m$ where the $x_i$'s are the endpoints of the edges that meet at the vertex. In words, the eigenfunctions are required to be
continuous at the vertex and the sum of their (outward) normal derivatives, $\partial_{n}f(x_i)$, at the vertex is zero. These are called Kirchoff or Neumann vertex conditions. The operator so defined
is self-adjoint \cite{berkolaiko2013introduction,kostrykin1999kirchhoff}. 

Our main interest is the spectrum of the standard Laplace operator $L$.  The spectrum if formed by all $ \lambda = k^2 $ such that the equation $ L \psi = \lambda \psi $
has a non-trivial solution. Every such function $ \psi $ is a solution to the differential equation $ - \psi'' (x) = \lambda \psi $ on every edge
and satisfies the vertex conditions at every vertex. The operator $L$ is nonnegative since its quadratic form is given by
$$ \langle L  u, u \rangle_{L_2(\Gamma)} = \int_\Gamma \vert u' (x) \vert^2 dx \geq 0. $$
The spectrum of the operator is discrete and is formed by a sequence of eigenvalues tending to $ + \infty. $

First of all we note that $ \lambda_0 = 0 $ is an eigenvalue with the eigenfunction $ \psi_0 (x) \equiv 1 $. It is easy to see that this eigenfunction is unique, provided $ \Gamma $ is connected.

Consider now any $ \lambda = k^2 > 0 .$
 The solution of the differential equation on each edge is given by a linear combination of $ e^{i k x}$ and $e^{-i k x}$. Imposing the vertex conditions on the eigenfunctions gives a linear equation system which has a solution if a certain $k$-dependent determinant, $D(k)$, is zero. One possible way is to consider these equations as a certain scattering process arriving to the following secular equation
 \begin{equation}
 D(k) := \det \left( S_e (k) S_v - I \right) = 0,
 \end{equation}
 where $ S_e (k) $ is the edge scattering matrix formed by $ 2 \times 2 $ matrices $$ \left(
 \begin{array}{cc}
 0 & e^{ik \ell_n} \\
 e^{ik \ell_n}  & 0
 \end{array} \right) ,$$
 $ \ell_n $ being the length of the corresponding edge $ E_n,$
  and $ S_v $ is the vertex scattering matrix formed by $ d_m \times d_m $ scattering matrices 
 $$ S^m = \left(
 \begin{array}{ccccc}
 -1 + 2/d_m & 2/d_m & 2/d_m & \dots & 2/d_m \\
 2/d_m & -1 + 2/d_m & 2/d_m & \dots & 2/d_m \\
 2/d_m & 2/d_m & -1 + 2/d_m & \dots & 2/d_m \\
 \vdots & \vdots & \vdots & \ddots & \vdots \\
  2/d_m & 2/d_m & 2/d_m & \dots & -1 + 2/d_m \\
 \end{array} \right) ,$$
 with $ d_m $ being the valency of the vertex $ V_m. $
 Observe that the matrices $ S_e $ and $ S_v $ have block-diagonal forms in different bases. For example these marices for the graph depicted in Fig. 1 are as follows:
 $$ S_e (k) =
 \left(
 \begin{array}{cccccccc}
 0 & e^{ikc1} & 0 & 0 & 0 & 0 & 0 & 0 \\
 e^{ikc1} & 0 & 0 & 0 & 0 & 0 & 0 & 0 \\
      0 & 0 & 0 & e^{ik\pi} & 0 & 0 & 0 & 0 \\
  0 & 0 & e^{ik \pi} & 0 & 0 & 0 & 0 & 0 \\
   0 & 0 & 0 & 0 & 0 & e^{ikc1} & 0 & 0 \\
    0 & 0 & 0 & 0 & e^{ikc1} & 0 & 0 & 0 \\
     0 & 0 & 0 & 0 & 0 & 0 & 0 & e^{ikc2} \\
      0 & 0 & 0 & 0 & 0 & 0 & e^{ikc2} & 0 \end{array} \right) ,$$
      $$  S_v = \left(
      \begin{array}{cccccccc}
     -1/3 & 0 & 0 & 0 & 0 & 2/3 & 2/3 & 0 \\
     0 & 0 & 1 & 0 & 0 & 0 & 0 & 0 \\
     0 & 1 & 0 & 0 & 0 & 0 & 0 & 0 \\
     0 & 0 & 0 & 0 & 1 & 0 & 0 & 0 \\
     0 & 0 & 0 & 1 & 0 & 0 & 0 & 0 \\
     2/3 & 0 & 0 & 0 & 0 & -1/3 & 2/3 & 0 \\
     2/3 & 0 & 0 & 0 & 0 & 2/3 & -1/3 & 0 \\
     0 & 0 & 0 & 0 & 0 & 0 & 0 & 1 \\
\end{array} \right) .$$ 
 A description of how to obtain this determinant is given by Gutkin and Smilansky \cite{gutkin2001can} and in Kurasov and Nowaczyk \cite{kurasov2005inverse} (with an example by Nowaczyk \cite{nowaczyk2005inverse}) and we will not repeat this description here. The most important property of the function $ D(k) $ is that it is an analytic function (trigonometric polynomial) and all its zeroes
 are real having orders equal to the multiplicities of the corresponding eigenvalues, maybe with the exception of the eigenvalue $ \lambda_0 = 0. $
 
In what follows we are going to call the spectrum of the standard Laplace operator on a metric graph $ \Gamma $ by the spectrum of $ \Gamma.$
One important observation is that the spectrum scales in an easy way as the lengths of all edges are scaled simultaneously, hence it is natural to normalise the graphs so that they have total length equal to one, which we do in this paper.

 \section{Spectrum driven evolution}

 There is a wide literature devoted to different growth models for graphs. Most of these models are of a probabilistic or random nature. Some models, probably motivated by applications,
 use local rules \cite{vazquez2003growing}. One such example is crystal growth which is motivated by ice crystals and snow flakes. Our growth process is highly non-local, since 
 the spectra of quantum graphs are very sensitive to the overall graph structure and cannot be determined by a subgraph of the graph.  
 The spectrum of a graph is one of its most important global characteristic encoding information about its topology and geometry.  
 Our goals are: 
 
\begin{itemize}
\item Find graphs having particular eigenvalues or conditions on the eigenvalues. This is a hard inverse problem.
\item Get an intuitive feeling for graphs having particular conditions on the eigenvalues, possibly leading to new conjectures.
\end{itemize}

 Below we describe the most elementary variant of spectrum driven evolution and discuss some generalisations. 
 In order to determine such evolution, let us choose any {\bf target spectrum} $
  \sigma^* = \{ \mu_0, \mu_1, \mu_2, \dots, \mu_N \} \subset \mathbb R $ with $ \mu_0 = 0 $.
  The role of the initial condition is played by any finite compact metric graph $ \Gamma_0. $
  
  The evolution is determined using the following inductive procedure
 \begin{enumerate}
 \item {\bf Initiatialisation} 
 \newline We start by choosing an initial graph $ \Gamma_0 $. From this initial parent graph we will generate a sequence of graphs $ \Gamma_n $ using the
 following
  \item {\bf Evolution rule}: 
  \newline Given a parent graph $ \Gamma_n $, we construct the next graph $ \Gamma_{n+1} $
  in the sequence using following the steps: 
  \begin{itemize}
  \item Consider all possible graphs obtained by attaching a new edge of length $ \ell = 1/N $, where $N$ is the number of edges, to one of
  the vertices either as a pendant edge or between any two existing vertices;
    \item Scale the new graphs to meet the normalisation requirement of total length one;
  \item Calculate the spectra of all such graphs and choose the one with the spectrum closest to the target spectrum $
  \sigma^* = \{ \mu_0, \mu_1, \mu_2, \dots, \mu_N \}. $ Denote the optimal graph by $ \Gamma_{n+1}. $ In the unlikely case that several graphs that have spectra equally close to the target spectrum, we simply chose any one of them.
  \end{itemize}
  \item Repeating the evolution procedure, this time with $ \Gamma_{n+1} $ as the parent graph, one obtains an infinite sequence of metric graphs.
 \end{enumerate}
  The distance between the actual spectrum  $ \sigma (\Gamma_{n}) = \{ \lambda_0,  \lambda_1, \dots \} $ and the target spectrum $ \sigma^* $ is the usual Euclidean distance
  \begin{equation}
  \left( \sum_{k=0}^N \vert \lambda_k- \mu_k \vert^2 \right)^{1/2} . 
  \end{equation}
  
  Our main goal is to study properties of such a sequence and its possible convergence. It is natural to consider the index $ n $ as a certain time-parameter
  and speak about time evolution of the graph. This evolution certainly depends on the evolution rule, but the role of the initial state (graph $ \Gamma_0 $)
  is also important. 
  
    The formulated growth model can be generalised in many ways and it is not our goal to examine all such possibilities:
  \begin{itemize}
  \item the class of allowed graphs may be constrained by, for example, considering only trees;
  \item the notion of the distance between the two spectra can be different - any distance can be used;
  \item instead of adding one edge, one may add an arbitrary {\it small} graph;
  \item one may allow to delete existing edges or to collapse cycles;
  \item any point on an edge may be considered as a degree two vertex, hence it is natural to allow to add edges between any two points 
  on a metric graph;
  \item instead of having a target spectrum, one may choose another rule like choosing the graph maximising (or minimising) the spectral gap 
  or the entropy. $ $ 
  \end{itemize}

Since our understanding of topological and geometric properties of the graphs is far from being satisfactory, the search for rules that determine interesting evolution requires checking new ideas and  examining different possibilities.
Our studies show that the dynamics is extremely sensitive to the target spectrum which makes it into a challenging problem to choose the
target spectrum so that the graphs in the sequence have interesting properties. By choosing the evolution rule in a special way
one may expect the generated graph sequence will converge to a graph having certain extremal properties.

%
%
%
%
%
%
%
%
%
%
%
%
%
%

\section{Program use}

Our program is written in the computer language Mathematica \cite{Mathematica}. Graphs are specified in two different ways, either directly by Mathematica graph objects or by their adjacency matrix. The edge lengths are stored in the elements of the matrix. Four edges can be given symbolic lengths, $c1$, $c2$, $c3$ and $c4$. If the graph has a numerical length, the program normalizes the graph to length one, where the length is the sum of the length of the edges. The program obtains the eigenvalues, $\lambda=k^2$, using three different methods. In the first method we simply plot the value of $D(k)$ as a function of $k$ and obtain the roots visually. $D(k)$ may contain $c1$ and $c2$ and it is possible to adjust these parameters interactively using sliders and get a feeling for how the roots depend on these parameters. Figure 2 presents $D(k)$ for the graph shown in Figure 1. $c1=\pi$ and $c2$ has been varied from zero to $\pi$ to $50\pi$. It can be difficult to know if  $D(k)=0$ at a point or if it just becomes close to zero. This problem can be seen in the left panel in Fig. 2, where $D(k) \leq 0$ for all $k$. A further problem is that even if one knows that $D(k)=0$, due to  $D(k)$ changing sign, it may be difficult to know the multiplicity of the root without further investigations. 

The program can also compute the zeroes of $D(k)$ numerically and sometimes analytically. Analytical solutions are naturally very hard, and often impossible to obtain. However, surprisingly difficult cases can be solved. For instance the graph in Figure 1 has analytical solutions when $c1=c2=\pi$. The roots are:
\begin{align}
&k=8 \pi n \nonumber \\ 
&k= -(8/3) \pi+ 8 \pi n \nonumber \\
&k=(8/3) \pi+ 8 \pi n \nonumber \\
&k= -4 i \log \left(\frac{1}{12} \left(-3-\sqrt{33}-i \sqrt{102-6
   \sqrt{33}}\right)\right) + 8 \pi n \\
&k=-4 i \log \left(\frac{1}{12} \left(-3-\sqrt{33}+i \sqrt{102-6
   \sqrt{33}}\right)\right) + 8 \pi n \nonumber \\ 
&k=-4 i \log \left(\frac{1}{12} \left(-3+\sqrt{33}-i \sqrt{6
   \left(17+\sqrt{33}\right)}\right)\right) + 8 \pi n \nonumber \\
&k=-4 i \log \left(\frac{1}{12} \left(-3+\sqrt{33}+i \sqrt{6
   \left(17+\sqrt{33}\right)}\right)\right) + 8 \pi n \nonumber
\end{align}

where $n\in Z$. Our program seems to find analytical solutions for all graphs where the edge-lengths are rationally dependent (that is, the ratio between any two edge-lengths can be written as $p/q$ with $p$ and $q$ integers). This depends on the capabilities of Mathematica but so far we have found no exception to this observation. If the graph is large, or if there are edges where $p/q$ is very small, then the calculation might have to be aborted due to excessive computing time. Analytic solutions are never found if the edges of the graph are not rationally dependent. Furthermore, Mathematica does not explicitly give the multiplicities of the roots although they can be found by inspection of $D(k)$ .

Despite these shortcomings we have found our program to be very useful to gain understanding of the spectra of different graphs. The raison d'\'etre for our program is mostly to introduce experimentation into the field of quantum graphs and possible insights should be proven rigorously.

The program also supports the spectral growth of graphs. The user specifies up to three k-values (which determine the eigenvalues via $\lambda=k^2$), the initial graph and defines how growth should proceed. For example the growth may only add pendant edges, add edges between any two vertices, delete edges and so on. The evolution of the k-values is shown graphically.
It is also possible to maximise $\lambda_2/\lambda_1$ as a goal. Any type of goal can be defined by the user and many starting graphs can be set graphically. Having set the rules, the program then evolves the graphs for a certain number of user-defined steps.

\section{Experiments}

We will now describe some experiments that we have done. 
\begin{itemize}
  \item In the first experiment, we chose the target spectrum $\sigma^*=(0, \pi^2, (2\pi)^2)$ and the goal was to minimise $\Vert \sigma(\Gamma_{1i})-\sigma^* \Vert$. $\sigma^*$ is the initial part of the spectrum of a path graph and unsurprisingly the growth of the graphs proceeds by adding an edge to the parent  graph, simply extending it. See Fig. 4.
  \item In the second experiment the goal was to maximise $\lambda_{1}$, also known as the spectral gap. We then find that the sequence consists of graphs with a large number of edges compared with the number of vertices and include all complete graphs with equal length edges (equilateral complete graphs). It has been proven that the spectral gap becomes arbitrarily high for equilateral complete graphs when the number of edges increases \cite{kennedy2015spectral}. See Fig. 5.
  \item In the third experiment the goal was to maximise $\lambda_{2}/\lambda_{1}$. We then find that the sequence consists of graphs (often complete graphs) which are joined by one edge, also known as dumbbell graphs. From this experiment we conjecture that $\lambda_{2}/\lambda_{1}$ can become arbitrarily large for this class of graphs. We tested this and could reach $\lambda_{2}/\lambda_{1}>64$ before reaching the limit of our laptop computer. In fact testing with dumbbell graphs where two general graphs are joined by one edge we also find that $\lambda_{2}/\lambda_{1}$ can become very large. See Fig. 6.
  \item In the fourth experiment the goal was to maximise $\lambda_{1}$. The growth was here alternating between adding an edge and deleting an edge. In this case we prefer to use the term evolution rather than growth of the graphs. We started with a grid graph and after a few steps we ended up with a complete graph. After this stage the graphs alternate as illustrated in the Fig. 7. If we have the different goal of minimising $\lambda_{1}$ instead of maximising $\lambda_{1}$ the final graph will be a path graph, as expected. 
   \item In the fifth experiments we chose a target spectrum of\\  $\sigma^*=(0, (2\pi)^2, (3\pi)^2, (3\pi)^2)$ and the goal was to minimise $\Vert \sigma(\Gamma_{1i})-\sigma^* \Vert$. The growth was restricted to trees. Fig. 8 shows the evolution of the k-values as well as the resulting graph after 12 steps. The evolution manages to find a tree that has eigenvalues close to the target values.
\end{itemize} 

\begin{figure}
\centering
\includegraphics[width=0.7\textwidth]{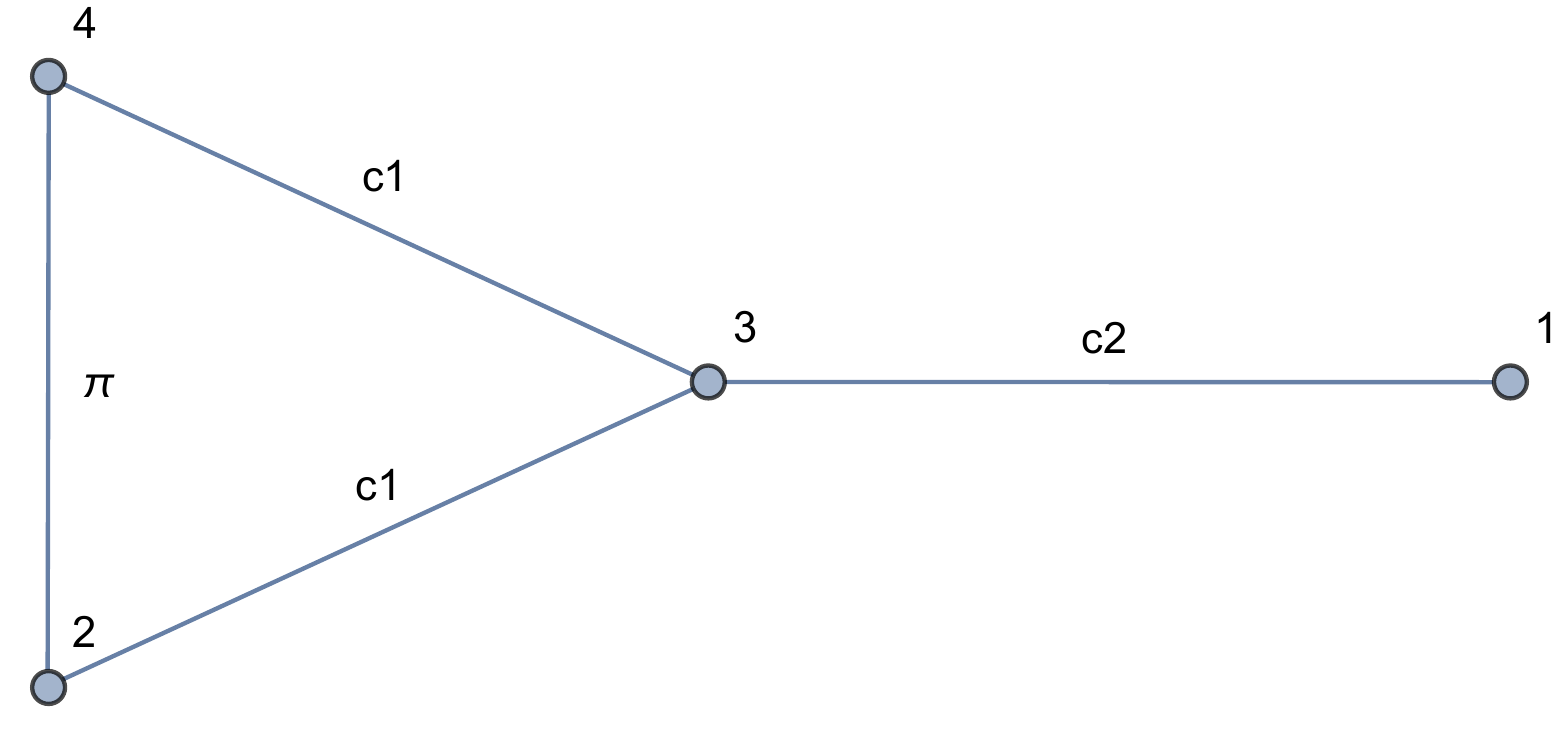}
\caption{\label{fig:graph}An example of a metric graph generated by our program from an adjacency matrix. The vertices are labeled by numbers and the edge lengths are indicated. $c1$ and $c2$ are  parameters.}
\end{figure}

\begin{figure}
\centering
\includegraphics[width=1.0\textwidth]{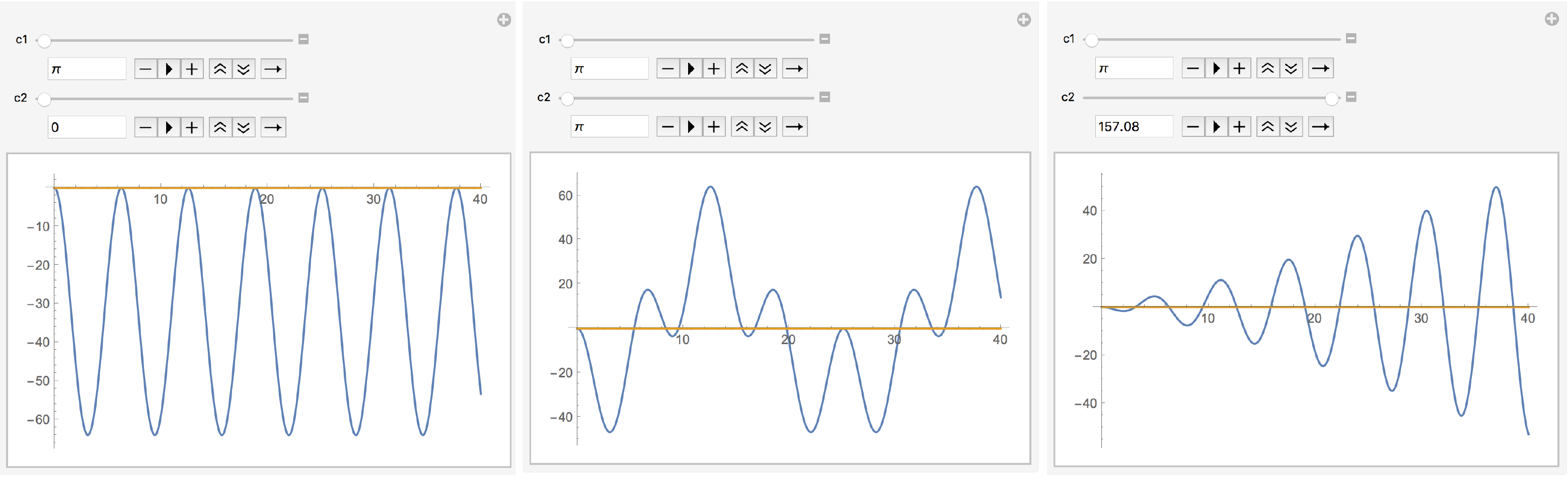}
\caption{\label{fig:all}Examples of plots of $D(k)$, for the graph shown in Fig. 1, where $c1$ and $c2$ can be changed using sliders. Left panel. Here $c1=\pi$ and $c2=0$ and we have a triangular graph (with length normalized to one). The zeros of $D(k)$ are then at $k=2n\pi$ which is reflected in the figure. Middle panel. Here $c1=c2=\pi$ and it is less obvious where the zeros of $D(k)$ should be located. Right panel. Here $c1=\pi$, $c2=50\pi$ and the graph is almost a path graph, which has zeros of $D(k)$ at $k=n\pi$.}
\end{figure}

\begin{figure}
\centering
\includegraphics[width=1.0\textwidth]{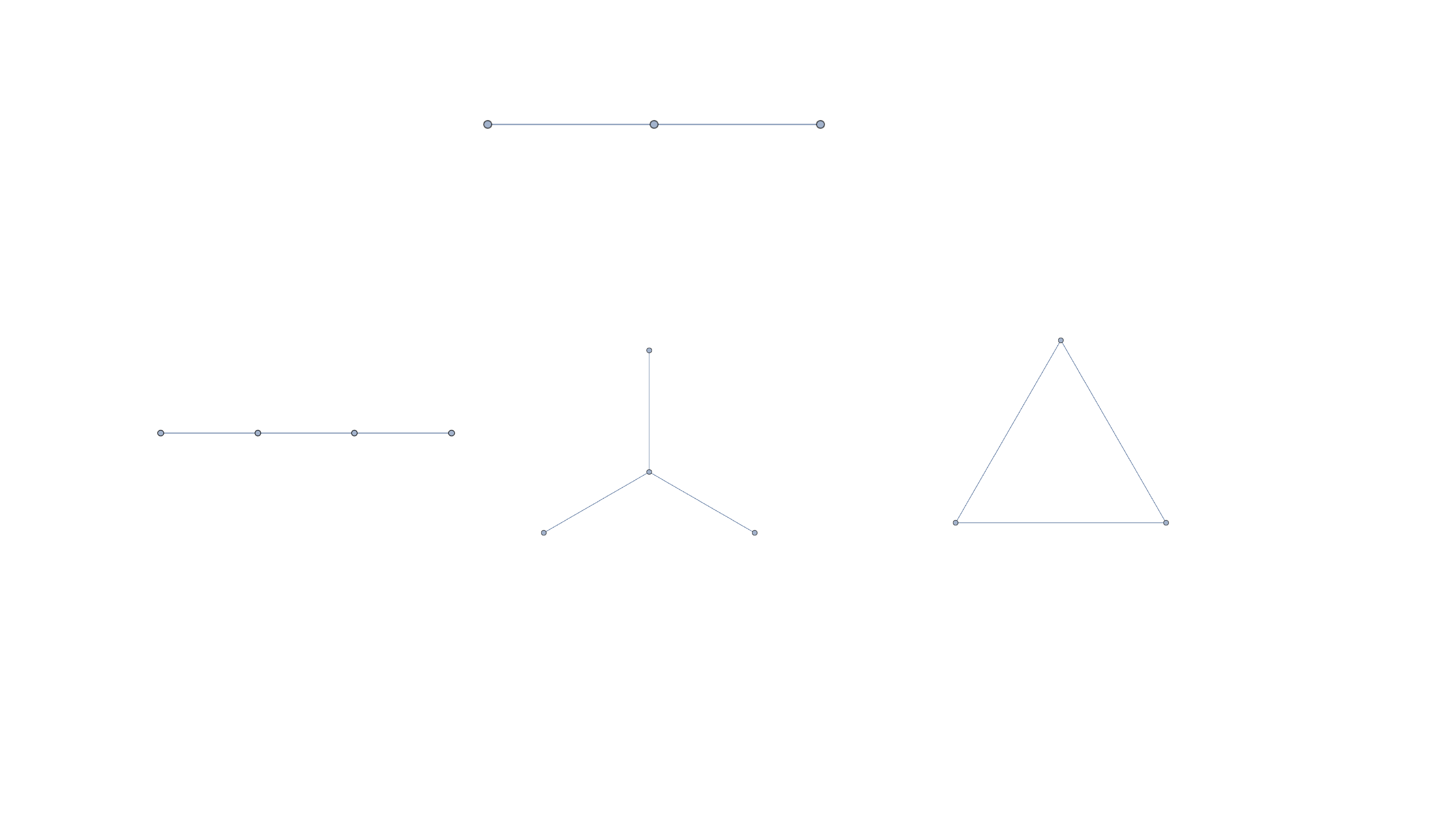}
\caption{\label{Parent-child:all}The top graph is the parent graph and has three vertices. From this graph we get three child graphs by adding an edge, shown in the bottom. Out of these child graphs we select the graph which has a spectrum which best fit some chosen criterion. The selected graph becomes a new parent graph and the process is repeated.}
\end{figure}

\begin{figure}
\centering
\includegraphics[width=1.0\textwidth]{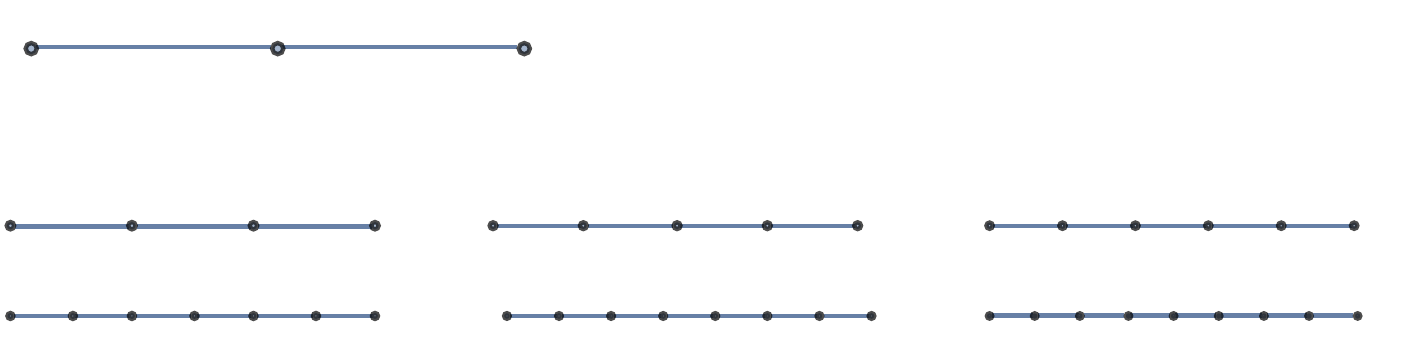}
\caption{\label{EdgeGrowth:all}The top graph is the parent graph and has three vertices. From this graph we get the sequence of graphs below. The goal was to have a spectrum with an initial segment which is close to $(0, \pi^2, (2 \pi)^2)$. This is satisfied by a path graph. The parent graph is simply extended when we evolve the graphs.}
\end{figure}

\begin{figure}
\centering
\includegraphics[width=1.0\textwidth]{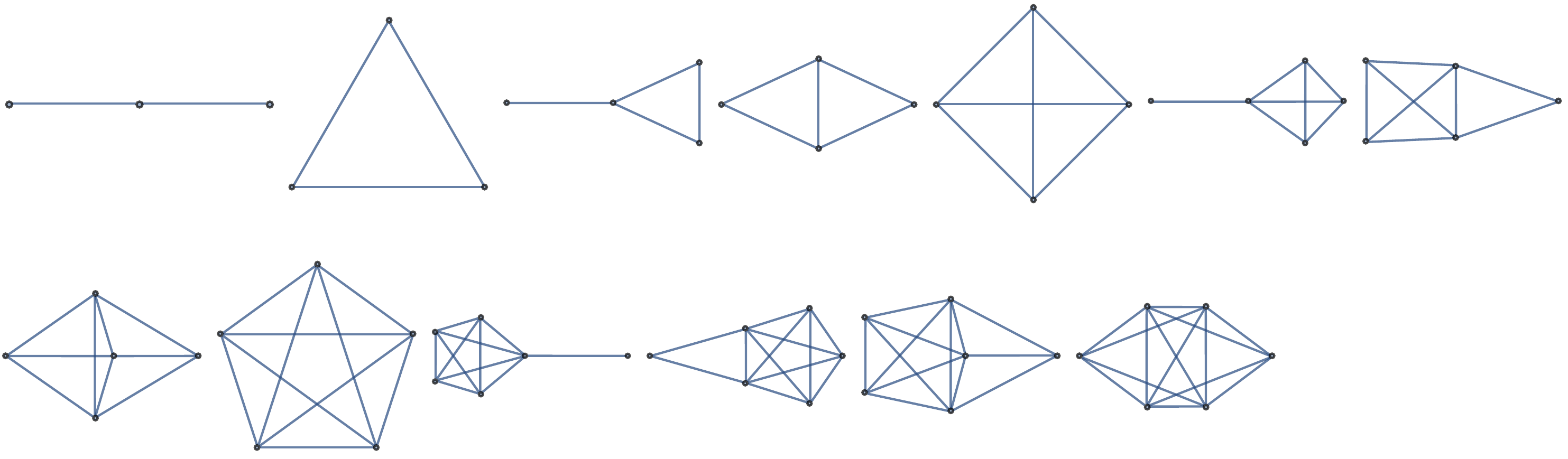}
\caption{\label{MaxSpectralGap:all}In this sequence of graphs the goal was to maximise the spectral gap, $\lambda_{1}$. The sequence contains many complete graphs.}
\end{figure}

\begin{figure}[h]
\centering
\includegraphics[width=1.0\textwidth]{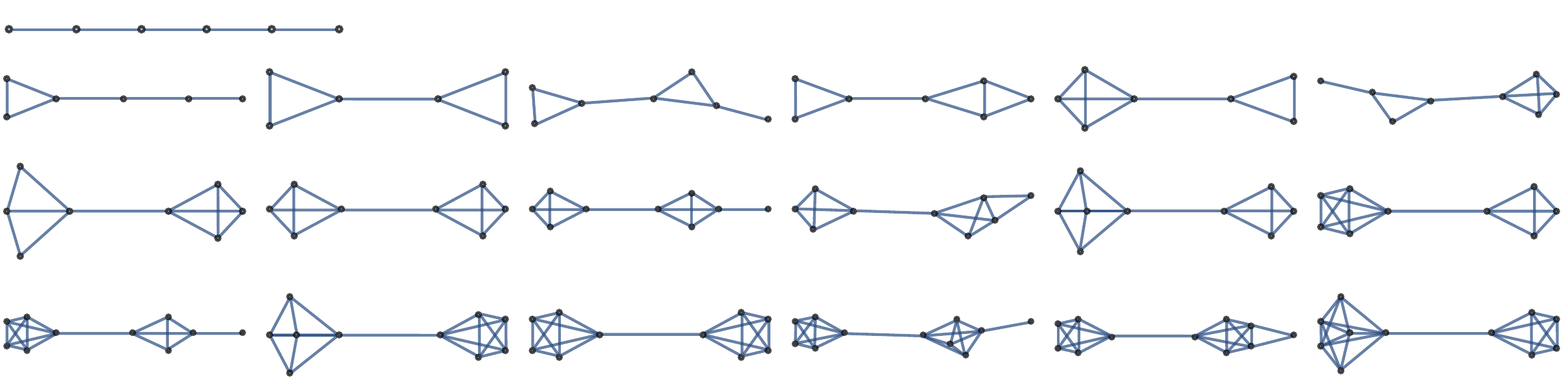}
\caption{\label{MaxRatio:all}In this sequence of graphs the goal was to maximise $\lambda_{2}/\lambda_{1}$. The graphs are two complete graphs joined by one edge. The parent graph is a path graph with six vertices. Not all graphs in the sequence have been plotted.}
\end{figure}

\begin{figure}
\centering
\includegraphics[width=1.0\textwidth]{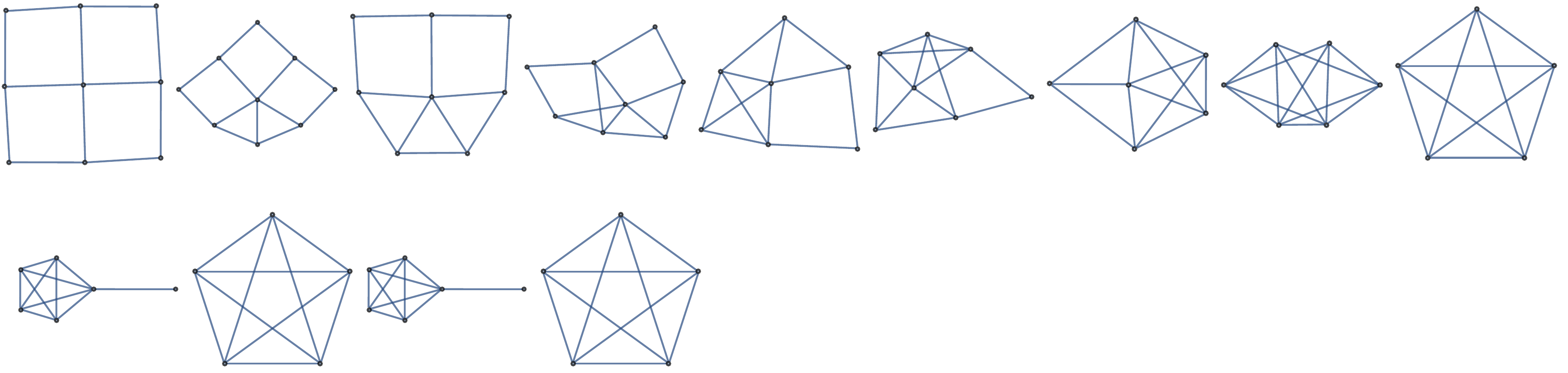}
\caption{\label{GridToComplete:all}In this sequence of graphs we added and then deleted an edge alternatingly, starting with a grid graph. The goal was to maximise the spectral gap, $\lambda_{1}$. The grid graph transforms to a complete graph after eight steps. The sequence then becomes cyclical, alternating between the complete graph and the complete graph with an added edge.}
\end{figure}

\begin{figure}
\centering
\includegraphics[width=0.4\textwidth]{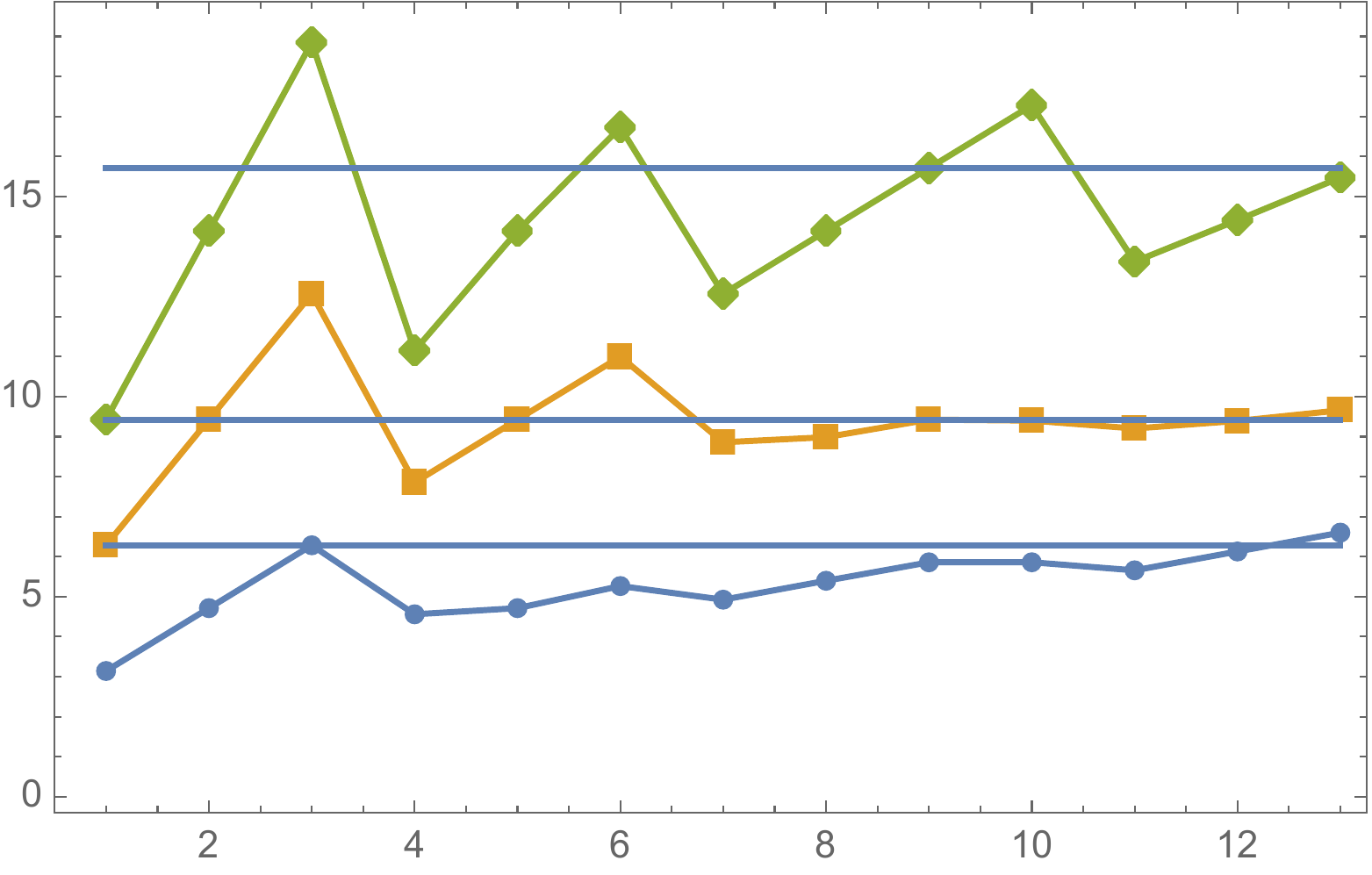}
\includegraphics[width=0.4\textwidth]{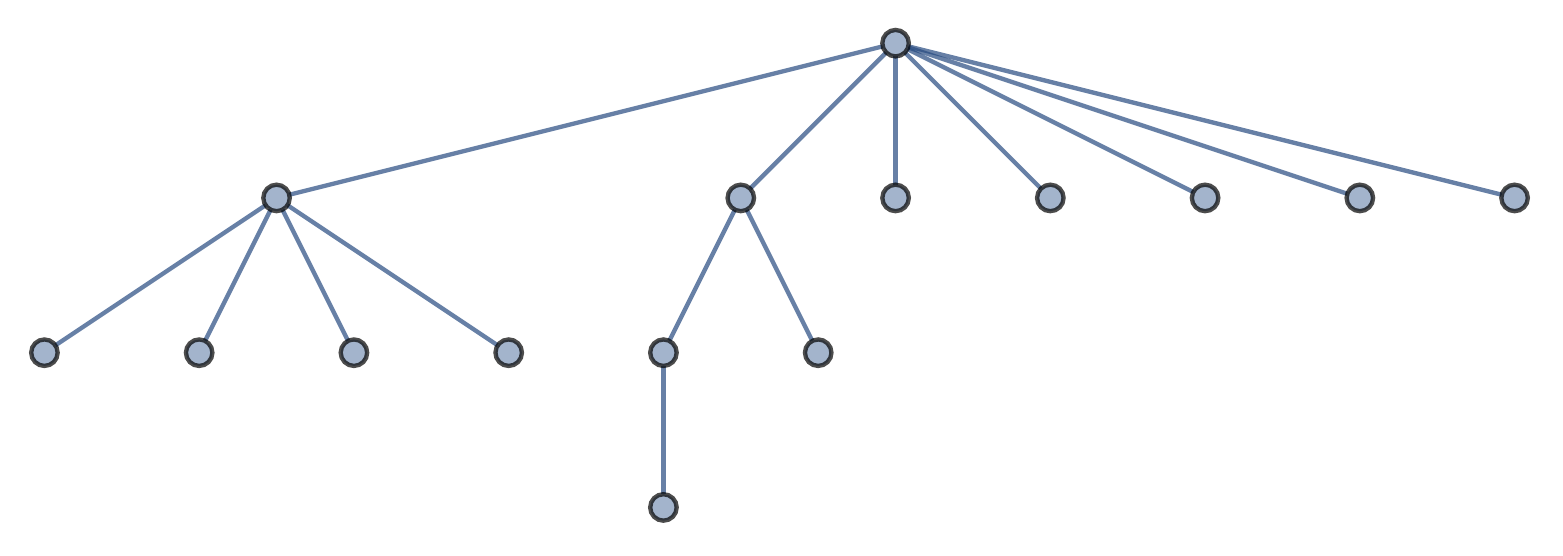}

\caption{\label{Eigenvalues:all}In this sequence the target spectrum was $(0, (2\pi)^2, (3\pi)^2, (3\pi)^2)$ and the growth was restricted to trees. We show the evolution of the k-values as well as the final graph in the sequence. The target k-values are indicated by horisontal lines. The spectral evolution manages to find a tree-graph having eigenvalues which are quite close to the target values.}
\end{figure}

\begin{figure}
\centering
\includegraphics[width=1.0\textwidth]{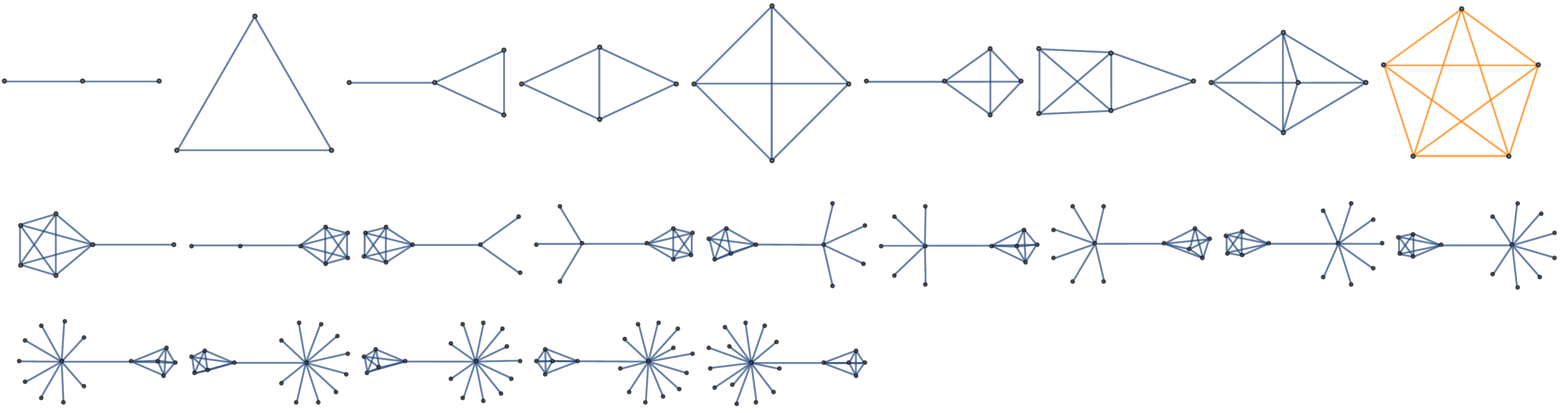}
\caption{\label{Program:all}In this sequence of graphs we used a program. The first goal was to increase $\lambda_{1}$ until $\lambda_{1}\geq (5\pi)^2$ by adding an edge. This resulted in a complete graph, yellow in the figure. The goal was then changed to have a spectrum close to $(0, (5\pi)^2,(15\pi)^2)$ by adding a vertex and an edge. The growth then proceeded as in the figure.}
\end{figure}

It is also possible to program the growth of the graphs. That is one can first specify a particular goal and how the evolution should proceed (e. g. add only pendant edges) and then switch goal and change how the evolution should proceed. In Fig. 9 we show an example of programmed growth of graphs.
By executing these type of programs one can steer the growth of the graphs to a certain goal such that subgoals are attained on the way. An interesting question is if it is possible to "help" the growth to attain a difficult goal by having suitable subgoals to attain first. We have found that some goals are very hard to obtain during our experiments. Unfortunately these types of programming experiments are very taxing for the computer and we could not execute too long programs.

\section{Conclusion and outlook}
We have found that studying the growth of graphs greatly helps our intuition about their spectral properties. One immediately learns that the eigenvalues of the graphs tend to increase with increasing complexity of the graphs. In fact the graph with the smallest $\lambda_1$ is the path graph \cite{nicaise1987spectre}, where $\lambda_1=\pi^2$, which most often is our starting graph. 

We also find that the growth of graphs can help to find patterns that can be made into conjectures and one example was given above, involving the first two non-trivial eigenvalues. Not much is known about the relation between graphs and more than one non-trivial eigenvalue. 
Future work likely will involve finding more conjectures, and possibly also proving them. The fact that our program sometimes can give all solutions in an analytical form is very helpful for such research. Concerning growth of graphs there is a huge number of questions to be answered. Is it for instance possible to have something analogous to phase transitions in the graphs? That is, can a small change in the goal completely change the topology of the graph?

We release our program as open-source \cite{Pistol2016}.  

\section{Acknowledgement}
This project was supported by the Swedish Research Council (VR) and NanoLund. 
 
\newpage
\bibliographystyle{unsrt}
\bibliography{bibliography.bib}

\end{document}